# COVARIANCE REGULARIZATION BY THRESHOLDING

By Peter J. Bickel[1] and Elizaveta Levina[2]

*University of California, Berkeley and University of Michigan*

This paper considers regularizing a covariance matrix of $p$ variables estimated from $n$ observations, by hard thresholding. We show that the thresholded estimate is consistent in the operator norm as long as the true covariance matrix is sparse in a suitable sense, the variables are Gaussian or sub-Gaussian, and $(\log p)/n \to 0$, and obtain explicit rates. The results are uniform over families of covariance matrices which satisfy a fairly natural notion of sparsity. We discuss an intuitive resampling scheme for threshold selection and prove a general cross-validation result that justifies this approach. We also compare thresholding to other covariance estimators in simulations and on an example from climate data.

**1. Introduction.** Estimation of covariance matrices is important in a number of areas of statistical analysis, including dimension reduction by principal component analysis (PCA), classification by linear or quadratic discriminant analysis (LDA and QDA), establishing independence and conditional independence relations in the context of graphical models, and setting confidence intervals on linear functions of the means of the components. In recent years, many application areas where these tools are used have been dealing with very high-dimensional datasets, and sample sizes can be very small relative to dimension. Examples include genetic data, brain imaging, spectroscopic imaging, climate data and many others.

It is well known by now that the empirical covariance matrix for samples of size $n$ from a $p$-variate Gaussian distribution, $\mathcal{N}_p(\mu, \Sigma_p)$, is not a good estimator of the population covariance if $p$ is large. Many results in random matrix theory illustrate this, from the classical Marčenko–Pastur law [29]

Received October 2007; revised February 2008.
[1]Supported by NSF Grant DMS-06-05236.
[2]Supported by NSF Grants DMS-05-05424 and DMS-08-05798 and NSA Grant MSPF-04Y-120.

to the more recent work of Johnstone and his students on the theory of the largest eigenvalues [12, 23, 30] and associated eigenvectors [24]. However, with the exception of a method for estimating the covariance spectrum [11], these probabilistic results do not offer alternatives to the sample covariance matrix.

Alternative estimators for large covariance matrices have therefore attracted a lot of attention recently. Two broad classes of covariance estimators have emerged: those that rely on a natural ordering among variables, and assume that variables far apart in the ordering are only weakly correlated, and those invariant to variable permutations. The first class includes regularizing the covariance matrix by banding or tapering [2, 3, 17], which we will discuss below. It also includes estimators based on regularizing the Cholesky factor of the inverse covariance matrix. These methods use the fact that the entries of the Cholesky factor have a regression interpretation, which allows application of regression regularization tools such as the lasso and ridge penalties [21], or the nested lasso penalty [28] specifically designed for the ordered variables situation. Banding the Cholesky factor has also been proposed [3, 34]. These estimators are appropriate for a number of applications with ordered data (time series, spectroscopy, climate data). For climate applications and other spatial data, since there is no total ordering on the plane, applying the Cholesky factor methodology is problematic; but as long as there is an appropriate metric on variable indexes (sometimes, simple geographical distance can be used), banding or tapering the covariance matrix can be applied.

However, there are many applications, for example, gene expression arrays, where there is no notion of distance between variables at all. These applications require estimators invariant under variable permutations. Shrinkage estimators are in this category and have been proposed early on [7, 20]. More recently, Ledoit and Wolf [26] proposed an estimator where the optimal amount of shrinkage is estimated from data. Shrinkage estimators shrink the overdispersed sample covariance eigenvalues, but they do not change the eigenvectors, which are also inconsistent [24], and do not result in sparse estimators. Several recent papers [5, 31, 35] construct a sparse permutation-invariant estimate of the *inverse* of the covariance matrix, also known as the concentration or precision matrix. Sparse concentration matrices are of interest in graphical models, since zero partial correlations imply a graph structure. The common approach of [5, 31, 35] is to add an $L_1$ (lasso) penalty on the entries of the concentration matrix to the normal likelihood, which results in shrinking some of the elements of the inverse to zero. In [31], it was shown that this method has a rate of convergence that is driven by $(\log p)/n$ and the sparsity of the truth. Computing this estimator is nontrivial for high dimensions and can be achieved either via a semidefinite programming algorithm [[5], [35]] or by using the Cholesky decomposition to reparametrize the concentration matrix [31], but all of these are computationally intensive. A faster algorithm that employs the lasso was proposed by Friedman, Hastie and Tibshirani [16]. This approach has also been extended to more general penalties like SCAD [15] by Lam and Fan [25] and Fan, Fan and Lv [14]. In specific applications, there have been other permutation-invariant approaches that use different notions of sparsity: Zou, Hastie and Tibshirani [36] apply the lasso penalty to loadings in PCA to achieve sparse representation; d'Aspremont et al. [6] compute sparse principal components by semidefinite programming; Johnstone and Lu [24] regularize PCA by moving to a sparse basis and thresholding; and Fan, Fan and Lv [13] impose sparsity on the covariance via a factor model, which is often appropriate in finance applications.

In this paper, we propose thresholding of the sample covariance matrix as a simple and permutation-invariant method of covariance regularization. This idea has been simultaneously and independently developed by El Karoui [10], who studied it under a special notion of sparsity called $\beta$-sparsity (see details in Section 2.4). Here we develop a natural permutation-invariant notion of sparsity which, though more specialized than El Karoui's, seems easier to analyze and parallels the treatment in [3] which defines a class of models where banding is appropriate. Bickel and Levina [3] showed that, uniformly over the class of approximately "bandable" matrices, the banded estimator is consistent in the operator norm (also known as the matrix 2-norm, or spectral norm) for Gaussian data as long as $(\log p)/n \to 0$.

Here we show consistency of the thresholded estimator in the operator norm, uniformly over the class of matrices that satisfy our notion of sparsity, as long as $(\log p)/n \to 0$, and obtain explicit rates of convergence. There are various arguments to show that convergence in the operator norm implies convergence of eigenvalues and eigenvectors [3, 10], so this norm is particularly appropriate for PCA applications. The rate we obtain is slightly worse than the rate of banding when the variables are ordered, but the difference is not sharp. This is expected, since in the situation when variables are ordered, banding takes advantage of the underlying true structure. Thresholding, on the other hand, is applicable to many more situations. In fact, our treatment is in many respects similar to the pioneering work on thresholding of Donoho and Johnstone [8] and the recent work of Johnstone and Silverman [22] and Abramovich et al. [1].

The rest of this paper is organized as follows. In Section 2 we introduce the thresholding estimator and our notion of sparsity, prove the convergence result and compare to results of El Karoui (Section 2.4) and to banding (Section 2.5). In Section 3, we discuss a cross-validation approach to threshold selection, which is novel in this context, and prove a cross-validation result of general interest. Section 4 gives simulations comparing several permutation-invariant estimators and banding. Section 5 gives an example of thresholding estimator applied to climate data and Section 6 gives a brief discussion. The Appendix contains more technical proofs.

**2. Asymptotic results for thresholding.** We start by setting up notation. We write $\lambda_{\max}(M) = \lambda_1(M) \geq \cdots \geq \lambda_p(M) = \lambda_{\min}(M)$ for the eigenvalues of a matrix $M$. Following the notation of [3], we define, for any $0 \leq r, s \leq \infty$ and a $p \times p$ matrix $M$,

$$\|M\|_{(r,s)} \equiv \sup\{\|M\mathbf{x}\|_s : \|\mathbf{x}\|_r = 1\}, \tag{1}$$

where $\|\mathbf{x}\|_r^r = \sum_{j=1}^p |x_j|^r$. In particular, we write $\|M\| = \|M\|_{(2,2)}$ for the operator norm, which for a symmetric matrix is given by

$$\|M\| = \max_{1 \leq j \leq p} |\lambda_j(M)|.$$

For symmetric matrices, we have (see, e.g., [18])

$$\|M\| \leq (\|M\|_{(1,1)} \|M\|_{(\infty,\infty)})^{1/2} = \|M\|_{(1,1)} = \max_j \sum_i |m_{ij}|. \tag{2}$$

We also use the Frobenius matrix norm,

$$\|M\|_F^2 = \sum_{i,j} m_{ij}^2 = \operatorname{tr}(MM^T).$$

We define the thresholding operator by

$$T_s(M) = [m_{ij} 1(|m_{ij}| \geq s)], \tag{3}$$

which we refer to as $M$ *thresholded at* $s$. Note that $T_s$ preserves symmetry and is invariant under permutations of variable labels, but does not necessarily preserve positive definiteness. However, if

$$\|T_s - T_0\| \leq \varepsilon \quad \text{and} \quad \lambda_{\min}(M) > \varepsilon, \tag{4}$$

then $T_s(M)$ is necessarily positive definite, since for all vectors $\mathbf{v}$ with $\|v\|_2 = 1$ we have $\mathbf{v}^T T_s M \mathbf{v} \geq \mathbf{v}^T M \mathbf{v} - \varepsilon \geq \lambda_{\min}(M) - \varepsilon > 0$.

2.1. *A uniformity class of covariance matrices.* Recall that the banding operator was defined in [3] as $B_k(M) = [m_{ij} 1(|i-j| \leq k)]$. The uniformity class of "approximately bandable" covariance matrices is defined by

$$\mathcal{U}(\varepsilon_0, \alpha, C) = \left\{ \Sigma : \max_j \sum_i \{|\sigma_{ij}| : |i-j| > k\} \leq C k^{-\alpha} \text{ for all } k > 0, \right.$$

$$\left. \text{and } 0 < \varepsilon_0 \leq \lambda_{\min}(\Sigma) \leq \lambda_{\max}(\Sigma) \leq 1/\varepsilon_0 \right\}. \tag{5}$$

Here we define the uniformity class of covariance matrices invariant under permutations by

$$\mathcal{U}_\tau(q, c_0(p), M) = \left\{ \Sigma : \sigma_{ii} \leq M, \sum_{j=1}^{p} |\sigma_{ij}|^q \leq c_0(p), \text{ for all } i \right\},$$

for $0 \leq q < 1$. Thus, if $q = 0$,

$$\mathcal{U}_\tau(0, c_0(p), M) = \left\{ \Sigma : \sigma_{ii} \leq M, \sum_{j=1}^{p} 1(\sigma_{ij} \neq 0) \leq c_0(p) \right\},$$

a class of sparse matrices. We will mainly write $c_0$ for $c_0(p)$ in the future. Note that

$$\lambda_{\max}(\Sigma) \leq \max_i \sum_j |\sigma_{ij}| \leq M^{1-q} c_0(p),$$

by the bound (2). Thus, if we define,

$$\mathcal{U}_\tau(q, c_0(p), M, \varepsilon_0) = \{\Sigma : \Sigma \in \mathcal{U}_\tau(q, c_0(p), M) \text{ and } \lambda_{\min}(\Sigma) \geq \varepsilon_0 > 0\},$$

we have a class analogous to (5).

Naturally, there is a class of covariance matrices that satisfies both banding and thresholding conditions. Define a subclass of $\mathcal{U}(\varepsilon_0, \alpha, C)$ by

$$\mathcal{V}(\varepsilon_0, \alpha, C) = \{\Sigma : |\sigma_{ij}| \leq C|i-j|^{-(\alpha+1)}, \text{ for all } i, j : |i-j| \geq 1,$$
$$\text{and } 0 < \varepsilon_0 \leq \lambda_{\min}(\Sigma) \leq \lambda_{\max}(\Sigma) \leq 1/\varepsilon_0\}.$$

for $\alpha > 0$. Evidently,

$$\mathcal{V}(\varepsilon, \alpha, C) \subset \mathcal{U}(\varepsilon_0, \alpha, C_1)$$

for $C_1 \leq C(1 + 1/\alpha)$.

On the other hand, $\Sigma \in \mathcal{V}(\varepsilon_0, \alpha, C)$ implies

$$\sum_j |\sigma_{ij}|^q \leq \varepsilon_0^{-q} + C \frac{(\alpha+1)q}{(\alpha+1)q - 1},$$

so that for a suitable choice of $c_0$ and $M$,

$$\mathcal{V}(\varepsilon_0, \alpha, C) \subset \mathcal{U}_\tau(q, c_0, M)$$

for $q > \frac{1}{\alpha+1}$.

2.2. *Main result.* Suppose we observe $n$ i.i.d. $p$-dimensional observations $\mathbf{X}_1, \ldots, \mathbf{X}_n$ distributed according to a distribution $F$, with $E\mathbf{X} = 0$ (without loss of generality), and $E(\mathbf{X}\mathbf{X}^T) = \Sigma$. We define the empirical (sample) covariance matrix by

$$\hat{\Sigma} = \frac{1}{n} \sum_{k=1}^n (\mathbf{X}_k - \bar{\mathbf{X}})(\mathbf{X}_k - \bar{\mathbf{X}})^T, \tag{6}$$

where $\bar{\mathbf{X}} = n^{-1} \sum_{k=1}^n \mathbf{X}_k$, and write $\hat{\Sigma} = [\hat{\sigma}_{ij}]$.

We have the following result which parallels the banding result (Theorem 1) of Bickel and Levina [3].

THEOREM 1. *Suppose $F$ is Gaussian. Then, uniformly on $\mathcal{U}_\tau(q, c_0(p), M)$, for sufficiently large $M'$, if*

$$t_n = M' \sqrt{\frac{\log p}{n}} \tag{7}$$

*and $\frac{\log p}{n} = o(1)$, then*

$$\|T_{t_n}(\hat{\Sigma}) - \Sigma\| = O_P\left(c_0(p) \left(\frac{\log p}{n}\right)^{(1-q)/2}\right)$$

*and uniformly on $\mathcal{U}_\tau(q, c_0(p), M, \varepsilon_0)$,*

$$\|(T_{t_n}(\hat{\Sigma}))^{-1} - \Sigma^{-1}\| = O_P\left(c_0(p) \left(\frac{\log p}{n}\right)^{(1-q)/2}\right).$$

PROOF. Recall that, without loss of generality, we assumed $E\mathbf{X} = \mathbf{0}$. Begin with the decomposition,

$$\hat{\Sigma} = \hat{\Sigma}^0 - \bar{\mathbf{X}}\bar{\mathbf{X}}^T, \tag{8}$$

where

$$\hat{\Sigma}^0 \equiv [\hat{\sigma}_{ij}^0] = \frac{1}{n} \sum_{k=1}^n \mathbf{X}_k \mathbf{X}_k^T.$$

Note that, by (8),

$$\max_{i,j} |\hat{\sigma}_{ij} - \sigma_{ij}| \leq \max_{i,j} |\hat{\sigma}_{ij}^0 - \sigma_{ij}| + \max_{i,j} |\bar{X}_i \bar{X}_j|. \tag{9}$$

By a result of Saulis and Statulevičius [32] adapted for this context in Lemma 3 of [3], and $\sigma_{ii} \leq M$ for all $i$,

$$P\left[\max_{i,j} |\hat{\sigma}_{ij}^0 - \sigma_{ij}| \geq t\right] \leq p^2 C_1 e^{-C_2 n t^2}, \tag{10}$$

COVARIANCE THRESHOLDING 7for $|t| < \delta$, where $C_1$, $C_2$ and $\delta$ are constants depending only on $M$. In particular, (10) holds if $t = o(1)$.

For the second term in (9), we have, by the union sum inequality, the Gaussian tail inequality and $\sigma_{ii} \leq M$ for all $i$,

$$\tag{11} P\left[\max_i |\bar{X}_i|^2 \geq t\right] \leq pC_3 e^{-C_4 nt}.$$

Combining (10) and (11), we see that if $\frac{\log p}{n} \to 0$ and $t = t_n$ is given by (7), then for $M'$ sufficiently large,

$$\tag{12} \max_{i,j} |\hat{\sigma}_{ij} - \sigma_{ij}| = O_P\left(\sqrt{\frac{\log p}{n}}\right).$$

We now recap an argument of Donoho and Johnstone [8]. Bound

$$\|T_t(\hat{\Sigma}) - \Sigma\| \leq \|T_t(\Sigma) - \Sigma\| + \|T_t(\hat{\Sigma}) - T_t(\Sigma)\|.$$

The first term above is bounded by

$$\tag{13} \max_i \sum_{j=1}^p |\sigma_{ij}| 1(|\sigma_{ij}| \leq t) \leq t^{1-q} c_0(p).$$

On the other hand,

$$\tag{14} \begin{aligned} \|T_t(\hat{\Sigma}) &- T_t(\Sigma)\| \\ &\leq \max_i \sum_{j=1}^p |\hat{\sigma}_{ij}| 1(|\hat{\sigma}_{ij}| \geq t, |\sigma_{ij}| < t) \\ &\quad + \max_i \sum_{j=1}^p |\sigma_{ij}| 1(|\hat{\sigma}_{ij}| < t, |\sigma_{ij}| \geq t) \\ &\quad + \max_i \sum_{j=1}^p |\hat{\sigma}_{ij} - \sigma_{ij}| 1(|\hat{\sigma}_{ij}| \geq t, |\sigma_{ij}| \geq t) \\ &= \mathrm{I} + \mathrm{II} + \mathrm{III}. \end{aligned}$$

Using (12), we have

$$\mathrm{III} \leq \max_{i,j} |\hat{\sigma}_{ij} - \sigma_{ij}| \max_i \sum_{j=1}^p |\sigma_{ij}|^q t^{-q} = O_P\left(c_0(p) t^{-q} \sqrt{\frac{\log p}{n}}\right).$$

To bound term I, write

$$\tag{15} \begin{aligned} \mathrm{I} &\leq \max_i \sum_{j=1}^p |\hat{\sigma}_{ij} - \sigma_{ij}| 1(|\hat{\sigma}_{ij}| \geq t, |\sigma_{ij}| < t) + \max_i \sum_{j=1}^p |\sigma_{ij}| 1(|\sigma_{ij}| < t) \\ &\leq \mathrm{IV} + \mathrm{V}. \end{aligned}$$

By (13),

(16) $$\text{V} \leq t^{1-q} c_0(p).$$

Now take $\gamma \in (0,1)$. Then,

$$\text{IV} \leq \max_i \sum_{j=1}^p |\hat{\sigma}_{ij} - \sigma_{ij}| 1(|\hat{\sigma}_{ij}| \geq t, |\sigma_{ij}| \leq \gamma t)$$

(17) $$+ \max_i \sum_{j=1}^p |\hat{\sigma}_{ij} - \sigma_{ij}| 1(|\hat{\sigma}_{ij}| > t, \gamma t < |\sigma_{ij}| < t)$$

$$\leq \max_{i,j} |\hat{\sigma}_{ij} - \sigma_{ij}| \max_i N_i(1-\gamma) + c_0(p)(\gamma t)^{-q} \max_{i,j} |\hat{\sigma}_{ij} - \sigma_{ij}|,$$

where $N_i(a) \equiv \sum_{j=1}^p 1(|\hat{\sigma}_{ij} - \sigma_{ij}| > at)$. Note that, for some $\delta > 0$,

(18) $$P\left[\max_i N_i(1-\gamma) > 0\right] = P\left[\max_{i,j} |\hat{\sigma}_{ij} - \sigma_{ij}| > (1-\gamma)t\right]$$
$$\leq p^2 e^{-n\delta(1-\gamma)^2 t^2},$$

if $t = o(1)$, uniformly on $\mathcal{U}$. By (18) and (16), and $0 < \gamma < 1$, if

(19) $$2\log p - n\delta t^2 \to -\infty,$$

then

(20) $$\text{IV} = O_P\left(c_0(p) t^{-q} \sqrt{\frac{\log p}{n}}\right)$$

and, by (9) and (13),

(21) $$\text{I} = O_P\left(c_0(p) t^{-q} \sqrt{\frac{\log p}{n}} + c_0(p) t^{1-q}\right).$$

For term II, we have

$$\text{II} \leq \max_i \sum_{j=1}^p [|\hat{\sigma}_{ij} - \sigma_{ij}| + |\hat{\sigma}_{ij}|] 1(|\hat{\sigma}_{ij}| < t, |\sigma_{ij}| \geq t)$$

(22) $$\leq \max_{i,j} |\hat{\sigma}_{ij} - \sigma_{ij}| \sum_{j=1}^p 1(|\sigma_{ij}| \geq t) + t \max_i \sum_{j=1}^p 1(|\sigma_{ij}| \geq t)$$

$$= O_P\left(c_0(p) t^{-q} \sqrt{\frac{\log p}{n}} + c_0(p) t^{1-q}\right).$$

Combining (21) and (22) and choosing $t$ as in (7) establishes the first claim of the theorem. The second claim follows since

$$\|[T_{t_n}(\hat{\Sigma})]^{-1} - \Sigma^{-1}\| = \Omega_P(\|T_{t_n}(\hat{\Sigma}) - \Sigma\|)$$

uniformly on $\mathcal{U}_\tau(q, c_0(p), M, \varepsilon_0)$, where $A = \Omega_P(B)$ means $A = O_P(B)$ and $B = O_P(A)$. $\square$

THEOREM 2. *Suppose $F$ is Gaussian. Then, uniformly on $\mathcal{U}_\tau(q, c_0(p), M)$, if $t = M'\sqrt{\frac{\log p}{n}}$ and $M'$ is sufficiently large,*

$$
(23) \qquad \frac{1}{p}\|T_t(\hat{\Sigma}) - \Sigma\|_F^2 = O_P\left(c_0(p)\left(\frac{\log p}{n}\right)^{1-q/2}\right).
$$

An analogous result holds for the inverse on $\mathcal{U}_\tau(q, c_0(p), M, \varepsilon_0)$. The proof of Theorem 2 is similar to the proof of Theorem 1 and can be found in the Appendix.

2.3. *The non-Gaussian case.* We consider two cases here. If, for some $\eta > 0$,

$$Ee^{tX_{ij}^2} \leq K < \infty \qquad \text{for all } |t| \leq \eta_j, \text{ for all } i, j,$$

then the proof goes through verbatim, since result (10) still holds. The bound on $\max_i |\bar{X}_i|^2$ will always be at least the squared rate of $\max_{i,j} |\hat{\sigma}_{ij} - \sigma_{ij}|$, hence we do not need normality for (11).

In the second case, if we have, for some $\gamma > 0$,

$$E|X_{ij}|^{2(1+\gamma)} \leq K \qquad \text{for all } i, j,$$

then by Markov's inequality

$$
(24) \qquad P[|\hat{\sigma}_{ij} - \sigma_{ij}| \geq t] \leq KC(\gamma)\frac{n^{-(1+\gamma)/2}}{t^{1+\gamma}}.
$$

Thus the bound (10) becomes

$$p^2 KC(\gamma)\frac{n^{-(1+\gamma)/2}}{t^{1+\gamma}}$$

and hence,

$$\max_{i,j}|\hat{\sigma}_{ij}^0 - \sigma_{ij}| = O_P\left(\frac{p^{2/(1+\gamma)}}{n^{1/2}}\right).$$

Therefore, taking $t_n = M\frac{p^{2/(1+\gamma)}}{n^{1/2}}$, we find that

$$
(25) \qquad \|T_{t_n}(\hat{\Sigma}) - \Sigma\| = O_P\left(c_0(p)\left(\frac{p^{2/(1+\gamma)}}{n^{1/2}}\right)^{1-q}\right),
$$

which is, we expect, minimax though this needs to be checked.

2.4. *Comparison to thresholding results of El Karoui.* El Karoui [10] shows as a special case that if:

(i) $E|X_j|^r < \infty$ for all $r$, $1 \leq j \leq p$,
(ii) $\sigma_{jj} \leq M < \infty$ for all $j$,
(iii) if $\sigma_{ij} \neq 0$, $|\sigma_{ij}| > Cn^{-\alpha_0}$, $0 < \alpha_0 < \frac{1}{2} - \delta_0 < \frac{1}{2}$,
(iv) $\Sigma$ is $\beta$-sparse, $\beta = \frac{1}{2} - \eta$, $\eta > 0$,
(v) $\frac{p}{n} \to c \in (0, \infty)$,

then, if $t_n = Cn^{-\alpha}$, $\alpha = \frac{1}{2} - \delta_0 > \alpha_0$,

$$\|T_{t_n}(\hat{\Sigma}) - \Sigma\| \xrightarrow{\text{a.s.}} 0. \qquad (26)$$

El Karoui's notion of $\beta$-sparsity is such that our case $q = 0$ is $\beta$-sparse with $c_0(p) = Kp^\beta$. Our results yield a rate of

$$O_P\left(\frac{p^{\beta + 2/(1+\gamma)}}{n^{1/2}}\right)$$

for $\gamma$ arbitrarily large. Since $\beta < \frac{1}{2}$ by assumption and $p \asymp n$ we see that our result implies (26) under (i), (ii), (iv), (v) and our notion of sparsity. Thus, our result is stronger than his in the all moments case, again under our stronger notion of sparsity. El Karoui's full result, in fact, couples a maximal value of $r$ in (i) with the largest possible value of $\beta$. Unfortunately, this coupling involves (iii) which we do not require. Nevertheless, his result implies the corresponding consistency results of ours, if (iii) is ignored, when only existence of a finite set of moments is assumed. However, according to El Karoui (personal communication), (iii) is not needed for (26) in the case when our sparsity condition holds.

2.5. *Comparison to banding results of Bickel and Levina.* Comparison is readily possible on $\mathcal{V}(\varepsilon_0, \alpha, C)$. By Theorem 1 of [3] the best rate achievable using banding is

$$O_P\left(\left(\frac{\log p}{n}\right)^{\alpha/(2(\alpha+1))}\right).$$

On the other hand, by our Theorem 1, thresholding yields

$$O_P\left(\left(\frac{\log p}{n}\right)^{(1-q)/2}\right)$$

for $q > \frac{1}{\alpha+1}$. Comparing exponents, we see that banding is slightly better in the situation where labels are meaningful, since we must have

$$1 - q < \frac{\alpha}{\alpha+1}.$$

However, since $1 - q$ can be arbitrarily close to $\frac{\alpha}{\alpha+1}$ the difference is not sharp. Not surprisingly, as $\alpha \to \infty$, the genuinely sparse case, the bounds both approach $(\frac{\log p}{n})^{1/2}$.

**3. Choice of threshold.** The question of threshold selection seems to be hard to answer analytically. In fact, the $\hat{\sigma}_{ij}$ have variances which depend on the distribution of $(X_i, X_j)$ through higher-order moments so it may in fact make sense to threshold differentially. We conjecture that this would not make much difference if we assume second and fourth moments bounded above and below. Ignoring this issue, we propose a cross-validation method analogous to the one used by Bickel and Levina [3] but made using the Frobenius metric which enables us to partly analyze it.

3.1. *Method.* Split the sample randomly into two pieces of size $n_1$ and $n_2$ where a choice to be "justified" theoretically is $n_1 = n(1 - \frac{1}{\log n})$, $n_2 = \frac{n}{\log n}$ and repeat this $N$ times. Let $\hat{\Sigma}_{1,\nu}$, $\hat{\Sigma}_{2,\nu}$ be the empirical covariance matrices based on the $n_1$ and $n_2$ observations, respectively, from the $\nu$th split. Form

$$\hat{R}(s) = \frac{1}{N} \sum_{\nu=1}^{N} \|T_s(\hat{\Sigma}_{1,\nu}) - \hat{\Sigma}_{2,\nu}\|_F^2 \tag{27}$$

and choose $\hat{s}$ to minimize $\hat{R}(s)$ (in practice for $s \geq \varepsilon_n \to 0$, $\varepsilon_n \asymp \frac{\log p}{n}$). We will show that, under the conditions of Theorem 2,

$$\frac{1}{p}\|T_{\hat{s}}(\hat{\Sigma}) - \Sigma\|_F^2 = O_P\left[\left(\frac{\log p}{n}\right)^{1-q/2} c_0(p)\right], \tag{28}$$

uniformly on $\mathcal{U}_\tau(q, c_0(p), M)$ for $q > 0$. Claim (28) is weaker than the desired

$$\|T_{\hat{s}}(\hat{\Sigma}) - \Sigma\|_{(2,2)}^2 = O_P\left[\left(\frac{\log p}{n}\right)^{1-q} c_0(p)\right], \tag{29}$$

in terms of the norm, though the left-hand side of (28), the average of a set of eigenvalues, can be viewed as a reasonable proxy for the operator norm, the maximum of the same set of eigenvalues.

We begin with two essential technical results of independent interest.

3.2. *An inequality.* We note an inequality derivable from a classic one of Pinelis—see [33], for instance.

PROPOSITION 1. *Let $\mathbf{U}_1, \ldots, \mathbf{U}_n$ be i.i.d. p-variate vectors with $E|\mathbf{U}_1|^2 \leq K$, $E\mathbf{U}_1 = \mathbf{0}$. Let $\mathbf{v}_1, \ldots, \mathbf{v}_J$ be fixed p-variate vectors of length 1. Define for $\mathbf{x} \in R^p$*

$$\|\mathbf{x}\|_\mathbf{v} = \max_{1 \leq j \leq J} |\mathbf{v}_j^T \mathbf{x}|.$$

*Then,*

$$E\left\|\sum_{i=1}^{n} \mathbf{U}_i\right\|_\mathbf{v}^2 \leq Cn \log J E\|\mathbf{U}_1\|_\mathbf{v}^2, \tag{30}$$

*where $C$ is an absolute constant.*

PROOF. By symmetrization,
$$E\left\|\sum_{i=1}^n \mathbf{U}_i\right\|_{\mathbf{v}}^2 \leq 2E\max_j\left(\sum_{i=1}^n \varepsilon_i|(\mathbf{U}_i - \mathbf{U}_i')^T\mathbf{v}_j|\right)^2,$$

where $\mathbf{U}_i'$ are i.i.d. as $\mathbf{U}_i$ and independent of $\mathbf{U}_i$, and $\{\varepsilon_i\}$ are $\pm 1$ with probability $1/2$ and independent of $|(\mathbf{U}_i - \mathbf{U}_i')^T\mathbf{v}_j|$. Let
$$W_{ij} = |(\mathbf{U}_i - \mathbf{U}_i')^T\mathbf{v}_j|, \qquad a_{ij} = \frac{W_{ij}}{(\sum_{i=1}^n W_{ij}^2)^{1/2}}.$$

Then,
$$E\max_{1\leq j\leq J}\left(\sum_{i=1}^n \varepsilon_i W_{ij}\right)^2$$
$$\leq E\left\{E\left[\max_{1\leq j\leq J}\left(\sum_{i=1}^n a_{ij}\varepsilon_i\right)^2 \Big| \{W_{ij}: 1\leq i\leq n, 1\leq j\leq J\}\right] \max_{1\leq j\leq J}\sum_{i=1}^n W_{ij}^2\right\}$$
$$\leq Cn\log J E\max_{1\leq j\leq J}\sum_{i=1}^n W_{ij}^2,$$

by Pinelis' inequality [33]. Thus
$$E\max_{1\leq j\leq J}\left(\sum_{i=1}^n \varepsilon_i W_{ij}\right)^2 \leq Cn\log J E\max_{1\leq j\leq J}((\mathbf{U}_i - \mathbf{U}_i')^T\mathbf{v}_j)^2$$
$$\leq 2Cn\log J E\max_{1\leq j\leq J}(\mathbf{U}_1^T\mathbf{v}_j)^2. \qquad \square$$

3.3. *A general result on V-fold cross-validation.* We will prove our result for $N=1$ in (27). The nature of our argument in Theorem 3 is such that it is fairly easy to see that it applies to each term of the sum in (27) and thus holds not just for the "sample splitting" ($N=1$) procedure, but also for the general 2-fold cross-validation procedure that is given by (27), and in fact for more general $V$-fold cross-validation procedures.

Let $\mathbf{W}_1, \ldots, \mathbf{W}_{n+B}$ be i.i.d. $Q$-variate vectors with distribution $P$, with $E_P\mathbf{W} \equiv \boldsymbol{\mu}(P)$. Let $\hat{\boldsymbol{\mu}}_j$, $1\leq j\leq J$ be estimates of $\boldsymbol{\mu}$ based on $\mathbf{W}_1, \ldots, \mathbf{W}_n$. For convenience, in this section we write $|\mathbf{x}|^2 = \|\mathbf{x}\|_2^2 = \sum_{j=1}^Q x_j^2$ and $(\mathbf{x}, \mathbf{y}) = \mathbf{x}^T\mathbf{y}$. Let
$$L(\boldsymbol{\mu}, \mathbf{d}) = |\boldsymbol{\mu} - \mathbf{d}|^2.$$

The *oracle estimate* $\hat{\boldsymbol{\mu}}^o$ is defined by
$$\hat{\boldsymbol{\mu}}^o \equiv \arg\min_j |\boldsymbol{\mu}(P) - \hat{\boldsymbol{\mu}}_j|^2.$$

The *sample splitting* estimate $\hat{\boldsymbol{\mu}}^c$ is defined as follows. Let

$$\bar{\mathbf{W}}_B = \frac{1}{B} \sum_{j=1}^{B} \mathbf{W}_{n+j}.$$

Then,

$$\hat{\boldsymbol{\mu}}^c \equiv \arg\min_j |\bar{\mathbf{W}}_B - \hat{\boldsymbol{\mu}}_j|^2.$$

Here is our basic result which has in some form appeared in Gyorfi et al. ([19], Chapter 7, Theorem 7.1, page 101), Bickel, Ritov and Zakai [4] and Dudoit and van der Laan [9]. The major public proof in [19] appears to be in error and does not directly apply to our case so we give the proof of our statement for completeness.

THEOREM 3. *Suppose:*

(A1) $|\hat{\boldsymbol{\mu}}^o - \boldsymbol{\mu}(P)|^2 = \Omega_p(r_n);$
(A2) $E_P \max_{1 \leq j \leq J} |(\mathbf{v}_j, \mathbf{W}_1 - \boldsymbol{\mu})|^2 \leq C\rho(J)$ *for any set* $\mathbf{v}_1, \ldots, \mathbf{v}_J$ *of unit vectors in* $\mathbb{R}^Q$;
(A3) $\rho(J_n) \frac{\log J_n}{B_n} = o(r_n)$.

*Then,*

(31) $$|\hat{\boldsymbol{\mu}}^c - \boldsymbol{\mu}(P)|^2 = |\hat{\boldsymbol{\mu}}^o - \boldsymbol{\mu}(P)|^2 (1 + o_P(1)) = \Omega_P(r_n),$$

*where* $A = \Omega_P(B)$ *means that* $A = O_P(B)$ *and* $B = O_P(A)$.

We identify suitable $J_n$ and $B_n$ in our discussion of Theorem 4.

PROOF OF THEOREM 3. By definition, writing $\boldsymbol{\mu} \equiv \boldsymbol{\mu}(P)$,

(32) $$|\hat{\boldsymbol{\mu}}^c - \bar{\mathbf{W}}_B|^2 \leq |\hat{\boldsymbol{\mu}}^o - \bar{\mathbf{W}}_B|^2,$$

which is equivalent to

(33) $$2(\hat{\boldsymbol{\mu}}^c - \hat{\boldsymbol{\mu}}^0, \bar{\mathbf{W}}_B - \boldsymbol{\mu}) \geq |\hat{\boldsymbol{\mu}}^c - \boldsymbol{\mu}|^2 - |\hat{\boldsymbol{\mu}}^o - \boldsymbol{\mu}|^2.$$

But,

(34) $$|(\hat{\boldsymbol{\mu}}^c - \hat{\boldsymbol{\mu}}^o, \bar{\mathbf{W}}_B - \boldsymbol{\mu})| \leq |(\hat{\boldsymbol{\mu}}^c - \boldsymbol{\mu}, \bar{\mathbf{W}}_B - \boldsymbol{\mu})| + |(\hat{\boldsymbol{\mu}}^o - \boldsymbol{\mu}, \bar{\mathbf{W}}_B - \boldsymbol{\mu})|.$$

Now, let

$$\hat{\boldsymbol{\nu}}_j = \frac{\hat{\boldsymbol{\mu}}_j - \boldsymbol{\mu}}{|\hat{\boldsymbol{\mu}}_j - \boldsymbol{\mu}|}.$$

Then we have

(35) $$|(\hat{\boldsymbol{\mu}}^c - \boldsymbol{\mu}, \bar{\mathbf{W}}_B - \boldsymbol{\mu})| \leq |\hat{\boldsymbol{\mu}}^c - \boldsymbol{\mu}| \max_{1 \leq j \leq J} |(\hat{\boldsymbol{\nu}}_j, \bar{\mathbf{W}}_B - \boldsymbol{\mu})|,$$

and similarly for the other term. Now, by Proposition 1 and assumption (A2),

$$E \max_{1 \leq j \leq J} |(\hat{\boldsymbol{\nu}}_j, \bar{\mathbf{W}}_B - \boldsymbol{\mu})|^2 \leq C \frac{\log J}{B} \rho(J), \tag{36}$$

where $C$ is used generically. Therefore, after some algebra and Cauchy–Schwarz, by (32),

$$|\hat{\boldsymbol{\mu}}^c - \boldsymbol{\mu}|^2 \leq O_P\left(\frac{\log^{1/2} J}{B^{1/2}} \rho^{1/2}(J)\right)(|\hat{\boldsymbol{\mu}}^c - \boldsymbol{\mu}| + |\hat{\boldsymbol{\mu}}^o - \boldsymbol{\mu}|) + |\hat{\boldsymbol{\mu}}^o - \boldsymbol{\mu}|^2. \tag{37}$$

Letting $|\hat{\boldsymbol{\mu}}^c - \boldsymbol{\mu}|^2 = a_n$, we can rewrite (33) as

$$a_n \leq C \frac{\log^{1/2} J}{B^{1/2}} \rho^{1/2}(J)(a_n^{1/2} + r_n^{1/2}) + r_n, \tag{38}$$

with probability $1 - \varepsilon(C)$, with $\varepsilon(C) \to 0$ as $C \to \infty$. Using (iii),

$$a_n \leq a_n^{1/2} o_P(r_n^{1/2}) + r_n(1 + o_P(1)).$$

But by definition,

$$a_n \geq r_n$$

and hence,

$$a_n^{1/2} \leq o_P(r_n^{1/2}) + r_n^{1/2}(1 + o_P(1))$$

and the theorem follows. □

We proceed to show the relevance of Theorem 3 in our context. As we indicated, it is enough to consider $N = 1$, and for convenience write the observations as

$$\mathbf{X}_1, \ldots, \mathbf{X}_m, \ldots, \mathbf{X}_{m+B},$$

where $n = m + B$. Form $\hat{\Sigma}_1$ and $\hat{\Sigma}_2$, the sample covariances of $\mathbf{X}_1, \ldots, \mathbf{X}_m$ and $\mathbf{X}_{m+1}, \ldots, \mathbf{X}_{m+B}$, respectively, and the estimates $T_o(\hat{\Sigma})$, $T_{\hat{s}}(\hat{\Sigma})$ corresponding to the oracle and statistician. By Theorem 2, it is clear that for all $\Sigma \in \mathcal{U}_\tau$,

$$\|T_o(\hat{\Sigma}) - \Sigma\|_F^2 = O_P\left(c_0(p) p \left(\frac{\log p}{n}\right)^{1-q/2}\right) \tag{39}$$

and the same holds for $T_o(\hat{\Sigma}^0)$, the oracle estimate applied to the covariance matrix computed with known means.

Let the optimizing $\hat{s}$ and the oracle $s$ be, in fact, obtained by searching over a grid $\{j\sqrt{\frac{\log p}{n}} : 0 \leq j \leq J_n\}$. For selected $\Sigma \in \mathcal{U}_\tau$, $O_P$ in (39) can be turned into $\Omega_P$. To see this, consider, for example, $\sigma_{ii} = 1$, $\sigma_{ij} = \varepsilon |i-j|^{-(\alpha+1)}$ for $i \neq j$. For $\varepsilon > 0$ sufficiently small, this matrix is positive definite and we can take the right-hand side of (39) to be $r_n \equiv Kc_0(p)p(\frac{\log p}{n})^{1-q/2}$ for some $K$. For $t = M'\sqrt{\frac{\log p}{n}}$ and $M'$ sufficiently large, we know that (42) holds as an identity, which yields a contribution of at least $r_n$. The remaining terms in the risk for $T_t(\hat{\Sigma})$ can only increase it. The same holds for $\hat{\Sigma}^0$. For $\Sigma$ such as this with $\inf_t \|T_t(\hat{\Sigma}) - \Sigma\|_F^2 = \inf_t \|T_t(\hat{\Sigma}^0) - \Sigma\|_F^2 = r_n$, we can state the following theorem. For simplicity of notation, in what follows we assume $c_0(p) \equiv c_0 < \infty$. The general case follows by simply rescaling $\hat{\Sigma}$ by $c_0(p)$.

THEOREM 4. *Suppose $\mathbf{X}_i$ are Gaussian, $\Sigma \in \mathcal{U}_\tau(q, c_0(p), M)$, and $O_P = \Omega_P$ in (39). Then, if $B_n = n\varepsilon(n,p)$, $(\log J)^3 = o(n^{q/2} c_0(p)(\log p)^{1-q/2} \varepsilon(n,p))$, then*

$$(40) \qquad \|T_{\hat{s}}(\hat{\Sigma}) - \Sigma\|_F = \|T_o(\hat{\Sigma}) - \Sigma\|_F (1 + o_P(1)).$$

*Thus, $\sup\{\|T_{\hat{s}}(\hat{\Sigma}) - \Sigma\|_F^2 : \Sigma \in \mathcal{U}_\tau(q, c_0(p), M)\} = Kc_0(p)p(\frac{\log p}{n})^{1-q/2}$, which is the optimal rate for the oracle as well.*

The proof of Theorem 4, which consists of several lemmas that allow us to apply the general Theorem 3, is given in the Appendix.

NOTES. 1. Evidently, with $\varepsilon(n,p) \sim (\log n)^{-1}$ we can take $J \sim n^\kappa$, for any $\kappa < \infty$ if $q > 0$, and if $p \sim n^\delta$, even if $q = 0$.

2. Similar results can be obtained for banding.

3. The assumption of Gaussianity can be relaxed. By applying Corollary 4.10 from Ledoux [27], we can include distributions $F$ of $\mathbf{X} = A\boldsymbol{\varepsilon}$, where $\boldsymbol{\varepsilon} = (\varepsilon_1, \ldots, \varepsilon_p)^T$ and the $\varepsilon_j$ are i.i.d. $|\varepsilon_j| \leq c < \infty$ (thanks to N. El Karoui for pointing this out).

**4. Simulation results.** The simulation results we present focus on comparing banding, thresholding, and two more permutation-invariant estimators: the sample covariance and the shrinkage estimator of Ledoit and Wolf [26]. We consider the AR(1) population covariance model,

$$(41) \qquad \Sigma = [\sigma_{ij}] = [\rho^{|i-j|}]$$

with $\rho = 0.7$. The value of 0.7 was chosen so that the matrix is not very sparse (as would be the case with $\rho \leq 0.5$) but does have a fair number of very small entries (which would not be the case with $\rho$ close to 1). For banding, we show results for the variables in their "correct" order, and permuted at random. All other estimators are invariant to variable permutations, so their results are the same for both of these scenarios. We consider three values of $p = 30, 100, 200$ and the sample size is fixed at $n = 100$.

TABLE 1
*Averages and standard errors over 100 replications of performance measures for AR(1) with $\rho = 0.7$*

| $p$ | Sample | Ledoit–Wolf | Banding | Banding perm. | Thresholding |
|---|---|---|---|---|---|
| | | | Matrix 1-norm | | |
| 30 | 3.87(0.07) | 3.36(0.05) | 2.54(0.05) | 3.85(0.07) | 3.28(0.05) |
| 100 | 11.46(0.09) | 7.99(0.05) | 3.13(0.04) | 5.05(0.01) | 4.61(0.04) |
| 200 | 22.00(0.14) | 11.82(0.06) | 3.34(0.03) | 5.09(0.01) | 4.99(0.01) |
| | | | Operator norm | | |
| 30 | 1.95(0.04) | 1.69(0.03) | 1.38(0.03) | 1.92(0.04) | 1.90(0.04) |
| 100 | 4.16(0.05) | 3.06(0.02) | 1.68(0.02) | 4.63(0.003) | 3.15(0.03) |
| 200 | 6.68(0.06) | 3.80(0.01) | 1.80(0.02) | 4.67(0.002) | 3.64(0.02) |
| | | | Frobenius norm | | |
| 30 | 3.19(0.04) | 2.89(0.03) | 2.42(0.03) | 3.21(0.04) | 3.42(0.03) |
| 100 | 10.23(0.04) | 8.16(0.02) | 4.60(0.02) | 13.80(0.001) | 8.73(0.03) |
| 200 | 20.24(0.05) | 14.02(0.02) | 6.61(0.03) | 19.61(0.001) | 13.79(0.03) |
| | | Abs. difference between true and estimated largest eigenvalue | | | |
| 30 | 0.91(0.06) | 0.46(0.04) | 0.52(0.04) | 0.84(0.06) | 0.74(0.05) |
| 100 | 2.86(0.06) | 0.43(0.03) | 0.38(0.03) | 4.24(0.01) | 1.07(0.05) |
| 200 | 5.21(0.07) | 0.42(0.03) | 0.31(0.02) | 4.23(0.01) | 1.15(0.04) |
| | | Abs. cosine of the angle between true and estimated 1st PC | | | |
| 30 | 0.77(0.03) | 0.77(0.03) | 0.81(0.01) | 0.76(0.03) | 0.70(0.02) |
| 100 | 0.37(0.02) | 0.37(0.02) | 0.42(0.02) | 0.10(0.004) | 0.28(0.01) |
| 200 | 0.27(0.02) | 0.27(0.02) | 0.26(0.01) | 0.06(0.003) | 0.18(0.01) |

Table 1 shows average losses and standard deviations over 100 replications, as measured by three different matrix norms (matrix 1-norm which we denote $\|\cdot\|_{(1,1)}$, operator and Frobenius norms). We also report the absolute difference in the largest eigenvalue, $|\lambda_{\max}(\hat{\Sigma}) - \lambda_{\max}(\Sigma)|$, and the absolute value of the cosine of the angle between the estimated and true eigenvectors corresponding to the first eigenvalue. This assesses how accurate each of the estimators would be in estimating the first principal component.

TABLE 2
*Averages and standard errors over 100 replications of selected band width and threshold*

| $p$ | Banding $k$ | Banding perm. $k$ | Threshold $t$ |
|---|---|---|---|
| 30 | 4.36(0.07) | 24.57(0.14) | 0.33(0.004) |
| 100 | 4.27(0.05) | 0.00(0) | 0.49(0.002) |
| 200 | 4.22(0.04) | 0.00(0) | 0.55(0.001) |

The results in Table 1 show what one would expect: when banding is given the correct order of variables, it performs better than thresholding, since it is taking advantage of the underlying structure. When banding is given the variables in the wrong order, it performs poorly, often worse than the sample covariance matrix, and then thresholding is a much better choice. The Ledoit–Wolf estimator performs worse than thresholding by most measures, although it does well on estimating the largest eigenvalue. Note that the eigenvectors of the Ledoit–Wolf estimator are equal to the sample covariance eigenvectors.

Table 2 shows the band width selected by the cross-validation procedure on correct and permuted orderings, and the threshold selection. Note that banding in permuted order always selects a diagonal model for both $p = 100$ and $p = 200$, and keeps almost all the entries at $p = 30$, both of which result in bad estimators. The selected threshold increases with dimension, which is expected since in higher dimensions one would need to regularize more. The selected band width also goes down with dimension (the decrease on our range of $p$'s is not very large, but results over a wider range of dimensions in [3] show the same pattern more clearly).

Figure 1 shows scree plots of the true eigenvalues and means, 2.5% and 97.5% percentiles of the estimates over 100 replications for $p = 100$. Interestingly, the results show that the sample covariance is very bad at estimating the leading eigenvalues, but better than thresholding on the middle part of the spectrum. The leading eigenvalues, however, are more important in applications like PCA. The Ledoit–Wolf estimator does better on eigenvalues than on overall loss measures in Table 1. Banding in the correct order appears to do best on estimating the spectrum. For illustration purposes, scree plots from a single randomly selected realization are shown in Figure 2.

**5. Climate data example.** In this section, we illustrate the performance of the thresholded covariance estimator by applying it to climate data. The data are monthly mean temperatures recorded from January 1850 to June 2006; only the January data were used in the analysis below (157 observations). The region covered by the total of 2592 recording stations extends from $-177.5$ to $177.5$ degrees longitude, and from $-87.5$ to $87.5$ latitude. Not all the stations were in place for the entire period of 157 years; we do not impute the missing data in any way, but instead simply calculate spatial covariance from all the years available for any given pair of stations.

EOFs (empirical orthogonal functions) are frequently used in spatio-temporal statistics to represent patterns in spatial data. They are simply the principal components of the spatial covariance matrix, where observations over time are used as replications to calculate covariance between different spatial locations. EOFs are typically represented by spatial contour plots, which provide a visual illustration of which regions contribute the most to which principal components.

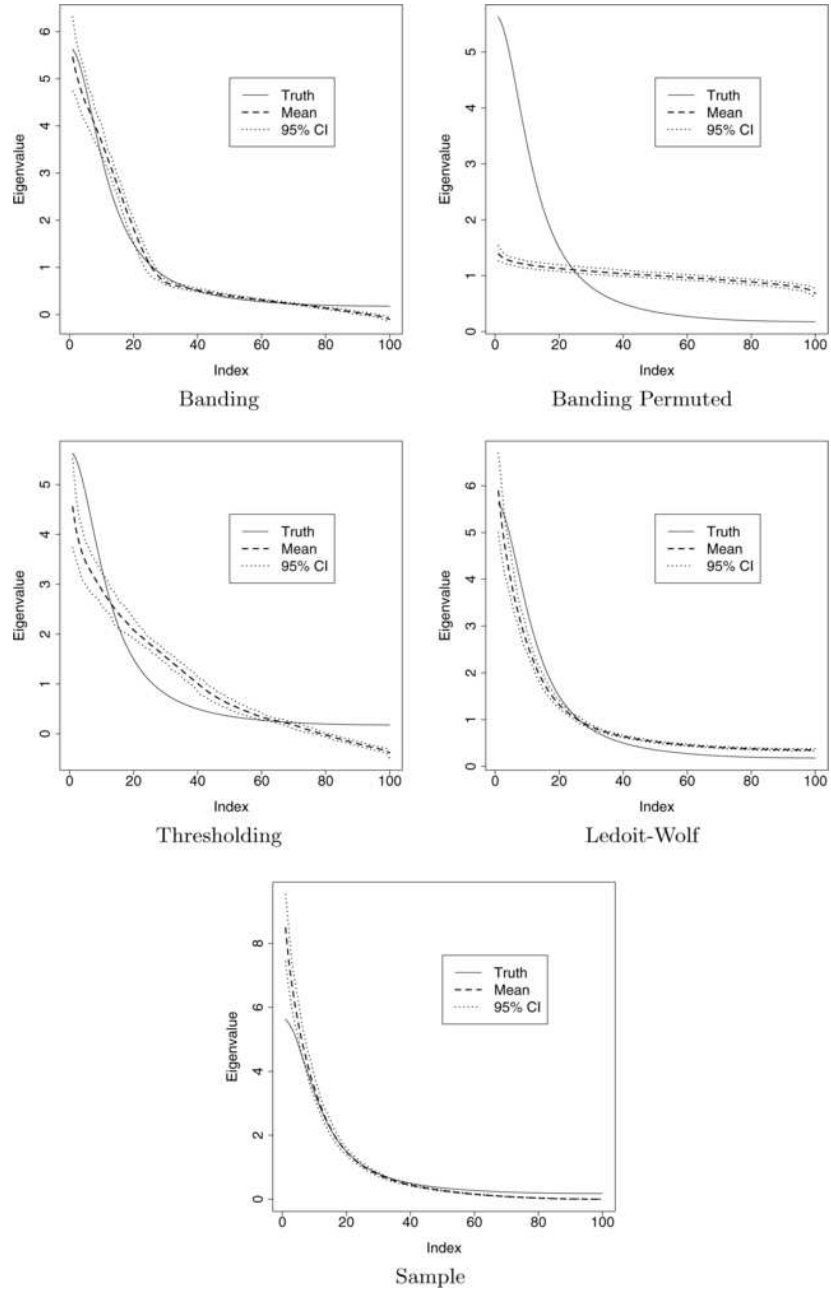

Fig. 1. *Scree plots: the mean estimated eigenvalues, their 2.5% and 97.5% percentiles, and the truth.*

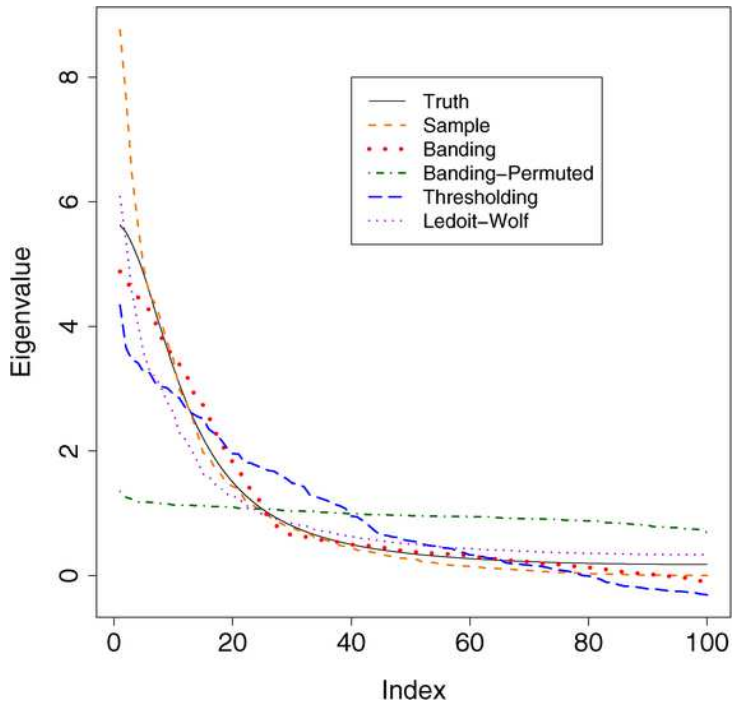

Fig. 2. *Scree plot of single realization.*

The plots in Figures 3 and 4 show the contour plots of the first four EOFs obtained, respectively, from the spatial sample covariance matrix (regular PCA), and from the thresholded spatial covariance matrix. We see that with thresholding, the first EOF essentially corresponds to Eurasia, and the second to North America, which the climate scientists agree should be separate. The regular PCA does not separate the continents. Ideal separation would be achieved if the estimator was block-diagonal (no nonzero correlations between North America and Eurasia). The thresholding estimator is in fact not block-diagonal, but does set enough correlations to zero to achieve effective separation of two continents in the EOFs. We also note that in this case the thresholding estimator does have some small negative eigenvalues, but they correspond to a negligible fraction of variance.

**6. Summary and discussion.** We have proposed and analyzed a regularization by thresholding approach to estimation of large covariance matrices. One of its biggest advantages is its simplicity—hard thresholding carries no computational burden, unlike many other methods for covariance regularization. A potential disadvantage is the loss of positive definiteness—but since we show that for a suitably sparse class of matrices the estimator is consistent as long as $(\log p)/n \to 0$, the estimator will be positive definite with probability tending to 1. We show consistency in the operator norm, which guarantees consistency for principal components, hence we expect that PCA will be one of the most important applications of the method.

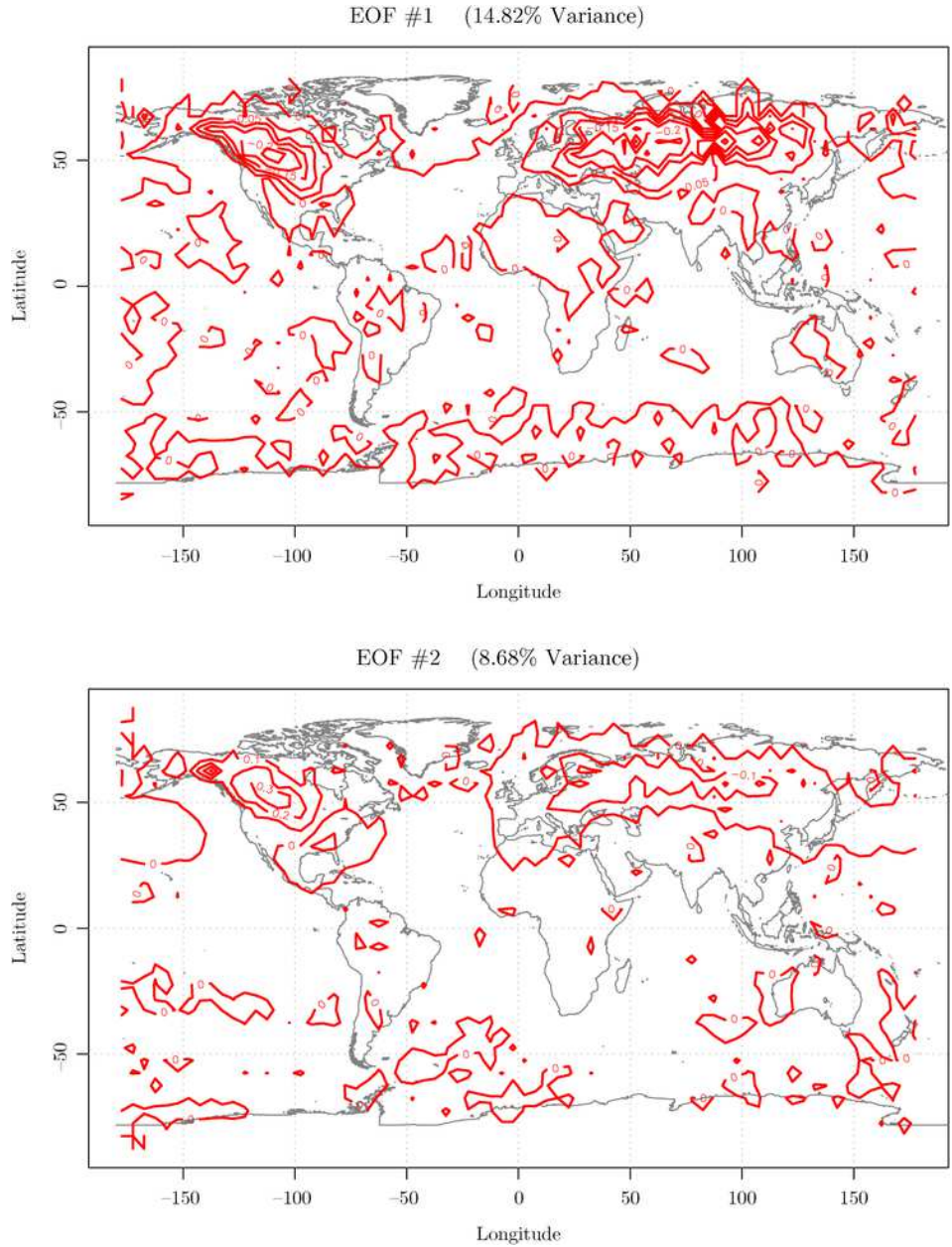

Fig. 3. *First two EOFs for the January temperature data obtained from regular PCA.*

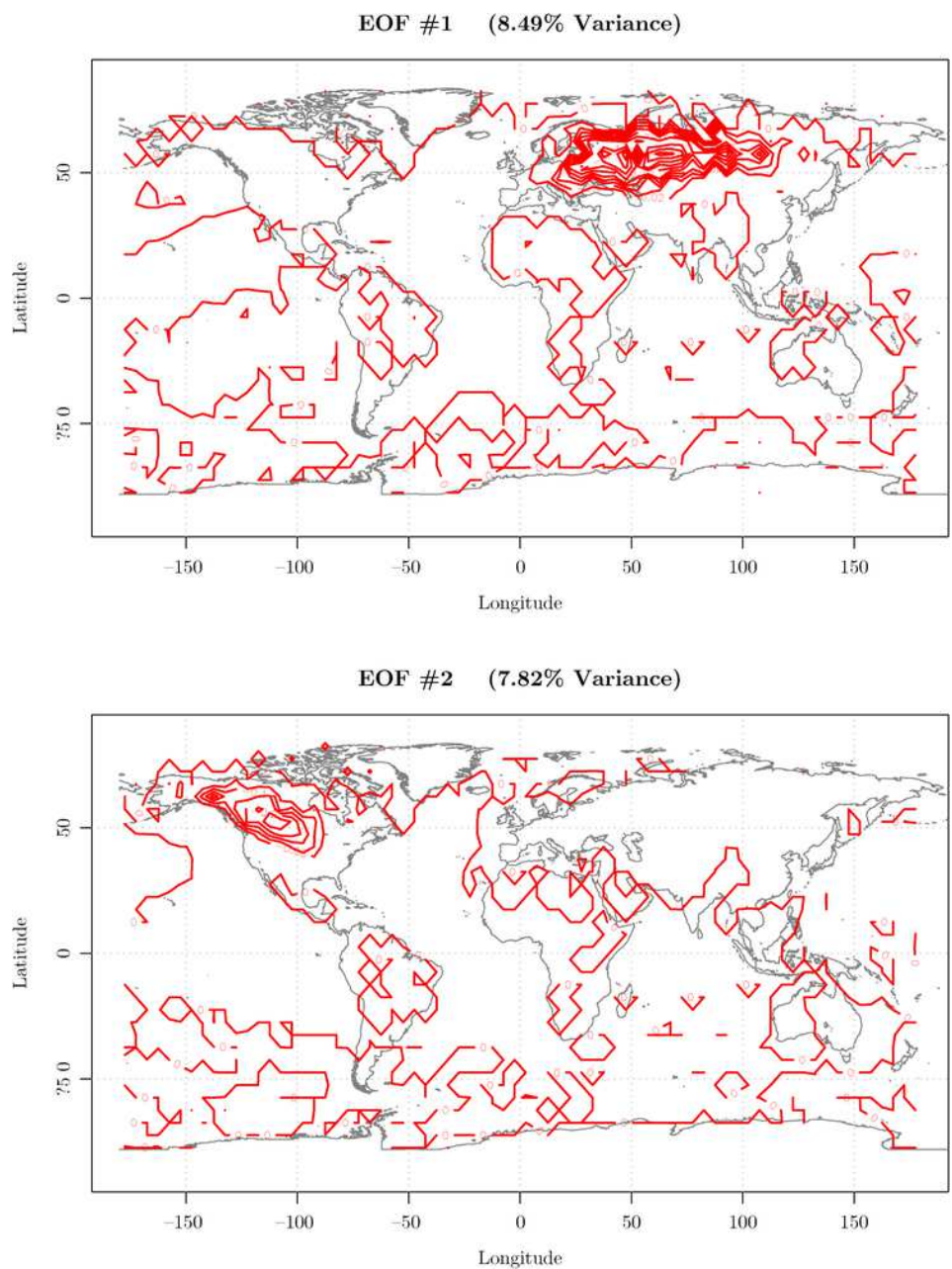

Fig. 4. *First two EOFs for the January temperature data obtained from PCA on the thresholded covariance matrix.*

We have also provided theoretical justification for the cross-validation approach to selecting the threshold. While it was formulated in the context of hard thresholding, the general result is much more widely applicable; in particular, it applies to other covariance estimation methods that depend on selecting the tuning parameter, such as [3, 21, 28, 31] and others.

## APPENDIX: ADDITIONAL PROOFS

PROOF OF THEOREM 2. The proof is essentially the same as for Theorem 1. We need to bound

$$\sum_{a,b}(\hat{\sigma}_{ab}1(|\hat{\sigma}_{ab}| \geq t) - \sigma_{ab})^2.$$

As before,

(42) $$\sum_{a,b} \sigma_{ab}^2 1(|\sigma_{ab}| < t) \leq t^{2-q} p c_0(p).$$

Similarly, for instance, the analogue of term III is

(43) $$\sum_{a,b}(\hat{\sigma}_{ab} - \sigma_{ab})^2 1(|\hat{\sigma}_{ab}| \geq t, |\sigma_{ab}| \geq t) \leq t^{-q} p c_0(p) \frac{\log p}{n}(1 + o_P(1))$$
$$= O_P\left(p c_0(p) \left(\frac{\log p}{n}\right)^{1-q/2}\right).$$

The only new type of term arises in the analogue of term I. Note that

$$\sum_{a,b} \hat{\sigma}_{ab}^2 1(|\hat{\sigma}_{ab}| \geq t, |\sigma_{ab}| < t)$$

$$\leq p c_0(p) t^{2-q} + \sum_{a,b} |\hat{\sigma}_{ab}^2 - \sigma_{ab}^2| 1(|\hat{\sigma}_{ab}| \geq t, |\sigma_{ab}| < t).$$

The second term is bounded by

(44) $$\Sigma_{a,b}|\hat{\sigma}_{ab}^2 - \sigma_{ab}^2|1(t \leq |\hat{\sigma}_{ab}| \leq (1+\varepsilon)t, |\sigma_{ab}| < t)$$
$$+ \Sigma_{a,b}|\hat{\sigma}_{ab}^2 - \sigma_{ab}^2|1(|\hat{\sigma}_{ab}| > (1+\varepsilon)t, |\sigma_{ab}| < t).$$

The second term in (44) is 0 with probability tending to 1. The first is bounded by

(45) $$\Sigma_{a,b}|\hat{\sigma}_{ab}^2 - \sigma_{ab}^2|1(t \leq |\hat{\sigma}_{ab}| \leq (1+\varepsilon)t, |\sigma_{ab}| \leq \gamma_1 t)$$
$$+ \Sigma_{a,b}|\hat{\sigma}_{ab}^2 - \sigma_{ab}^2|1(t \leq |\hat{\sigma}_{ab}| \leq (1+\varepsilon)t, \gamma_1 t < |\sigma_{ab}| < t),$$

where $\gamma_1 \in (0,1)$. Again, the first term in (45) is 0 with probability tending to 1 and the second is bounded by

$$(1+\varepsilon)t^2 c_0(p) p \gamma_1^{-q} t^{-q}$$

since

$$\max_{a,b} |\hat{\sigma}_{ab}^2 - \sigma_{ab}^2| \leq \max_{a,b} |\hat{\sigma}_{ab} + \sigma_{ab}| \max_{a,b} |\hat{\sigma}_{ab} - \sigma_{ab}|.$$

The theorem follows by putting (42), (43) and the other remainder terms together. It is clear from the argument that, by restricting the result from the class $\mathcal{U}_\tau(q, c_0(p), M)$ to properly chosen $\Sigma$'s, we can change $O_P$ into $\Omega_P$. □

We need the following lemmas to apply Theorem 3 to the special case of Theorem 4. Let

$$\mathbf{W}_i \equiv \{[X_a^{(i)} X_b^{(i)}], 1 \leq a, b \leq p\}$$

where $\mathbf{X}_i \equiv (X_1^{(i)}, \ldots, X_p^{(i)})^T$, so that $Q = p^2$.

LEMMA A.1. *Suppose $\Sigma \in \mathcal{U}_\tau(q, c_0, M)$, $V = [v_{ab}]$ is symmetric $p \times p$, $\|V\|_F = 1$ and $F$ is Gaussian. Then,*

$$\mathrm{Var}\left(\sum_{a,b} v_{ab} X_a^{(1)} X_b^{(1)}\right) \leq p C_1(q, c_0, M),$$

$$E\left(\sum_{a,b} v_{ab} \bar{X}_a \bar{X}_b\right)^2 \leq \frac{p}{B^2} C_2(q, c_0, M).$$

PROOF. By Wick's theorem, $E(X_a X_b X_c X_d) = \sigma_{ab}\sigma_{cd} + \sigma_{ac}\sigma_{bd} + \sigma_{ad}\sigma_{bc}$. Then for any $V = [v_{ab}]$ as above,

$$(46) \quad E\left(\sum_{a,b} v_{ab}(X_a^{(1)} X_b^{(1)} - \sigma_{ab})\right)^2 = \sum_{a,b,c,d} v_{ab} v_{cd}(\sigma_{ac}\sigma_{bd} + \sigma_{ad}\sigma_{bc}).$$

The two terms are equal because $v_{ab} = v_{ba}$. Consider the first term. Write

$$\left|\sum_{a,b,c,d} v_{ab} v_{cd} \sigma_{ac} \sigma_{bd}\right| = \left|\sum_{a,d}\left(\sum_c \sigma_{ac} v_{cd}\right)\left(\sum_b \sigma_{bd} v_{ab}\right)\right|$$

$$\leq \sum_{a,d}\left(\sum_c \sigma_{ac}^2\right)^{1/2}\left(\sum_c v_{cd}^2\right)^{1/2}\left(\sum_b \sigma_{bd}^2\right)^{1/2}\left(\sum_b v_{ab}^2\right)^{1/2}$$

$$\leq M^{2-q} c_0 \left(\sum_a \left(\sum_b v_{ab}^2\right)^{1/2}\right)^2 \leq M^{2-q} c_0 p,$$

where we used

$$\sum_c \sigma_{ac}^2 \leq \max_{a,c} |\sigma_{ac}|^{2-q} \sum_c |\sigma_{ac}|^q \leq M^{2-q} c_0.$$

Similarly,

$$E\left(\sum_{a,b} v_{ab} \bar{X}_a \bar{X}_b\right)^2$$

$$= \frac{1}{B^4} E\left(\sum_{a,b,c,d,j,k,l,m} v_{ab} v_{cd} X_a^{(j)} X_b^{(k)} X_c^{(l)} X_d^{(m)}\right)$$

$$= \frac{1}{B^3} \sum_{a,b,c,d} v_{ab} v_{cd} (\sigma_{ac}\sigma_{bd} + \sigma_{ad}\sigma_{bc}) + \frac{1}{B^2} \sum_{a,b,c,d} v_{ab} v_{cd} \sigma_{ab} \sigma_{cd}.$$

The first term is the same as (46). The second term is

$$\left(\sum_{a,b} v_{ab} \sigma_{ab}\right)^2 \leq \left(\sum_a \left(\sum_b \sigma_{ab}^2\right)^{1/2} \left(\sum_b v_{ab}^2\right)^{1/2}\right)^2 \leq M^{2-q} p c_0$$

and Lemma A.1 follows. $\square$

LEMMA A.2. *Suppose $\Sigma \in \mathcal{U}_\tau(q, c_0, M)$, $V$ is symmetric $p \times p$, $\|V\|_F = 1$, $\hat{\Sigma}_1^0 = [X_a X_b^T]$ and $F$ is Gaussian. Let $\Sigma = P\Lambda P^T$ be the eigendecomposition of $\Sigma$, and $|\gamma_1| \geq |\gamma_2| \geq \cdots$ be the eigenvalues of*

$$S = \Lambda^{1/2} P^T V P \Lambda^{1/2}.$$

*Then,*

$$(47) \qquad P\left[p^{-1/2}\left|\mathrm{tr}(V\hat{\Sigma}_1^0) - \sum_{j=1}^p \gamma_j\right| \geq t\right] \leq K e^{-\delta t(1+o(1))}$$

*for $t \to \infty$, $K = K(q, c_0, M) > 0$, $\delta = \delta(q, c_0, M) > 0$.*

PROOF. Let $\mathbf{Z} \sim \mathcal{N}_p(\mathbf{0}, I_{p \times p})$ and write

$$(48) \qquad \mathrm{tr}(V\hat{\Sigma}_1^0) = \mathrm{tr}(VP\Lambda^{1/2}\mathbf{Z}\mathbf{Z}^T\Lambda^{1/2}P^T) = \mathbf{Z}^T S \mathbf{Z} \sim \sum_{j=1}^p \gamma_j Z_j^2.$$

By Lemma A.1,

$$(49) \qquad \sum_{j=1}^p \gamma_j^2 = \tfrac{1}{2} \mathrm{Var}(\mathrm{tr}(V\hat{\Sigma}_1^0)) \leq p C(q, c_0, M).$$

Let $\tilde{\gamma}_j = \gamma_j p^{-1/2}$. In view of (48), to prove the right tail bound in (47) it is enough to show

$$P\left[\sum_{j=1}^{p} \tilde{\gamma}_j (Z_j^2 - 1) \geq t\right] \leq K e^{-\delta t(1+o(1))}$$

for $t \to \infty$ if

(50) $$\sum_{j=1}^{p} \tilde{\gamma}_j^2 \leq C(q, c_0, M).$$

Suppose first that all the $\gamma_j$ are positive and $\gamma_1 > \gamma_2$. By the general form of Markov's inequality,

(51) $$P\left[\sum_{j=1}^{p} \tilde{\gamma}_j (Z_j^2 - 1) \geq t\right] \leq \inf_{s} \left\{ \exp\left(-st - \sum_{j=1}^{p} \tilde{\gamma}_j s\right) \prod_{j=1}^{p} E e^{s \tilde{\gamma}_j Z_j^2} \right\}.$$

The log of the function on the right-hand side of (51) is

(52) $$-st + \sum_{j=1}^{p} [-\tilde{\gamma}_j s - \tfrac{1}{2} \log(1 - 2\tilde{\gamma}_j s)].$$

The minimizer satisfies

(53)
$$t = \sum_{j=1}^{p} [\tilde{\gamma}_j (1 - 2\tilde{\gamma}_j s)^{-1} - \tilde{\gamma}_j] = 2s \sum_{j=1}^{p} \tilde{\gamma}_j^2 (1 - 2\tilde{\gamma}_j s)^{-1}$$

$$= 2s \tilde{\gamma}_1^2 (1 - 2\tilde{\gamma}_1 s)^{-1} \left(1 + \sum_{j=2}^{p} \left(\frac{\tilde{\gamma}_j}{\tilde{\gamma}_1}\right)^2 (1 - 2\tilde{\gamma}_1 s)(1 - 2\tilde{\gamma}_j s)^{-1}\right).$$

Write $\omega \equiv 1 - 2\tilde{\gamma}_1 s$. Substituting into (53) and expanding:

(54) $$t = \frac{\tilde{\gamma}_1 (1 - \omega)}{\omega} \left(1 + \sum_{j=2}^{p} \tilde{\gamma}_j^2 \left(1 - \frac{\tilde{\gamma}_j}{\tilde{\gamma}_1}\right)^{-1} [\omega + O(\omega^2)]\right).$$

Now write $w \equiv \omega_0 + \Delta/t$, where $\omega_0 t = \tilde{\gamma}_1 (1 - \omega_0)$, and use $\omega_0 = (\tilde{\gamma}_1/t) + O(t^{-2})$. Substituting into (54) and solving for $\Delta$, we obtain after some computation,

$$\Delta = t^{-1} \sum_{j=2}^{p} \tilde{\gamma}_j^2 \left(1 - \frac{\tilde{\gamma}_j}{\tilde{\gamma}_1}\right)^{-1} + O(t^{-2}).$$

Note that $t\Delta$ is bounded by (50). Substituting back into (52) we finally obtain the bound

$$e^{(-t/(2\tilde{\gamma}_1))(1+o(1))}.$$

26    P. J. BICKEL AND E. LEVINAIf $\gamma_1$ has multiplicity $m > 1$, we can argue as above after pulling out $(1 - 2\tilde{\gamma}_1)^{-m/2}$ first. If the multiplicity of $\gamma_1$ is 1 but $\tilde{\gamma}_1 - \tilde{\gamma}_j \to 0$, the number of such terms is bounded above by (50) unless $\tilde{\gamma}_1$ itself tends to 0 and we can argue as if all such $\tilde{\gamma}_j$ are equal. Note also, that by (50) $m\tilde{\gamma}_1^2$ is bounded unless $\gamma_1 \to 0$. But in that case we can obtain a better bound than (47) by applying Theorem 3.1 of Saulis and Statulevičius [32].

Finally, if not all the $\gamma_j$ are positive, we break up

$$\sum_{j=1}^{p} \tilde{\gamma}_j(Z_j^2 - 1) = \sum_{j=1}^{p} \tilde{\gamma}_j^+(Z_j^2 - 1) - \sum_{j=1}^{p} \tilde{\gamma}_j^-(Z_j^2 - 1),$$

and compute the bound for each summand. A similar but easier argument and a better bound hold for $P[\sum_{j=1}^{p} \tilde{\gamma}_j(Z_j^2 - 1) \leq -t]$. The lemma follows. □

LEMMA A.3. *Under the conditions of Lemma A.2 for $\Sigma$ and each $V_j$,*

(55)  $$\rho(J) \leq C(q, c_0, M)(\log J)^2 p.$$

PROOF.    By Lemma A.2,

(56)
$$P\left[p^{-1/2} \max_{1 \leq j \leq J} \left| \operatorname{tr}(V_j \hat{\Sigma}_1^0) - \sum_{k=1}^{p} \gamma_{jk} \right| \geq t \right]$$
$$\leq \mathbf{1}(0 \leq t \leq x) + JKe^{-\delta t}\mathbf{1}(t > x)$$

by applying the union sum inequality for $t \geq x \to \infty$. Integrating, we get

$$p^{-1}E \max_{1 \leq j \leq J}\left(\operatorname{tr}(V_j \hat{\Sigma}_1^0) - \sum_{k=1}^{p} \gamma_{jk}\right)^2 \leq x^2 + JK\int_{x^2}^{\infty} e^{-\sqrt{t}\delta}\,dt$$

$$= x^2 + JK\int_{x}^{\infty} ve^{-v\delta}\,dv$$

$$= x^2 + JK(xe^{-x\delta} + \delta^{-1}e^{-x\delta}).$$

Minimizing over $x$, we get that as $J \to \infty$, the minimizer satisfies

$$x = A\log J(1 + o(1))$$

for $A(C, K) < \infty$. Then,

$$p^{-1}E \max_{1 \leq j \leq J}\left(\operatorname{tr}(V_j \hat{\Sigma}_1^0) - \sum_{k=1}^{p} \gamma_{jk}\right)^2 \leq C(\log J)^2$$

and Lemma A.3 follows.   □

A similar argument shows that under the conditions of Lemmas A.1 and A.2

$$(57) \quad E \max_{1 \leq j \leq J} \text{tr}\left[V_j\left(\bar{\mathbf{X}}\bar{\mathbf{X}}^T - \frac{1}{B}\Sigma\right)\right]^2 \leq B^{-2}C_3(q, c_0, M)p(\log J)^2.$$

PROOF OF THEOREM 4. We use Lemma A.3 to bound $\rho(J)$ and first obtain the equivalent of (40) for $\hat{\Sigma}^0$ by plugging in $\rho(J) \leq C(\log J)^2 p$. Take $r_n = Kp(\frac{\log p}{n})^{1-q/2}$, $B$, $J$ as in the statement of Theorem 4, and check that the conditions of Theorem 3 are satisfied when applied to $\{\mathbf{W}_i\}$ and the $J$ thresholding estimates that we optimize over for estimating $E(\mathbf{W}_1) = \Sigma$.

The argument for $\hat{\Sigma}$ requires us to return to (35) with $\bar{\mathbf{W}}_B - \boldsymbol{\mu}$ replaced by $\bar{\mathbf{W}}_B - \bar{\mathbf{X}}_B\bar{\mathbf{X}}_B^T$. By (57), (36) holds with this replacement on the left-hand side, and $\rho(J) \leq C(\log J)^2 pB^{-2}$. The analogue of Theorem 3 holds for this special case and Theorem 4 follows. $\square$

**Acknowledgments.** We thank Noureddine El Karoui, Sourav Chatterjee (UC Berkeley) and Iain Johnstone (Stanford) for helpful discussions; Adam Rothman (University of Michigan) for help with simulations; Yazhen Wang (University of Connecticut) for comments and a correction; Donghui Yan (UC Berkeley) and Serge Guillas (Georgia Tech University) for climate plots; and Gabi Hegerl (Duke University) for help in interpreting the climate plots.

Department of Statistics
University of California, Berkeley
Berkeley, California 94720-3860
USA
E-mail: bickel@stat.berkeley.edu

Department of Statistics
University of Michigan
Ann Arbor, Michigan 48109-1107
USA
E-mail: elevina@umich.edu